\documentclass[a4paper,12pt,twoside]{amsart}

\usepackage[cm]{fullpage}
\usepackage{mymacros,oitp}
\usepackage[pagebackref]{hyperref}
\title[Modular forms for (2,18)]{The ring of modular forms for the even unimodular lattice of signature (2,18)}
\author[A.~Nagano]{Atsuhira Nagano}
\address{
Faculty of Mathematics and Physics,
Institute of Science and Engineering,
Kanazawa University,
Kakuma, Kanazawa, Ishikawa,
920-1192, Japan}
\email{atsuhira.nagano@gmail.com}
\thanks{A.~N.~was partially supported by JSPS Kakenhi (18K13383) and MEXT LEADER}
\author[K.~Ueda]{Kazushi Ueda}
\address{
Graduate School of Mathematical Sciences,
The University of Tokyo,
3-8-1 Komaba,
Meguro-ku,
Tokyo,
153-8914,
Japan.}
\email{kazushi@ms.u-tokyo.ac.jp}
\thanks{K.~U.~was partially supported by JSPS Kakenhi (16H03930).}
\date{}

\begin{document}

\begin{abstract}
We show that the ring of modular forms with characters
for the even unimodular lattice of signature (2,18)
is obtained from the invariant ring of
$\mathrm{Sym}(\mathrm{Sym}^8(V) \oplus \mathrm{Sym}^{12}(V))$
with respect to the action of $\mathrm{SL}(V)$
by adding a Borcherds product of weight 132
with one relation of weight 264,
where $V$ is a 2-dimensional $\mathbb{C}$-vector space.
The proof is based on the study of
the moduli space of elliptic K3 surfaces with a section.
\end{abstract}

\maketitle

\section{Introduction}

Let $\bfU$
be the even unimodular hyperbolic lattice of rank 2.
A $\bfU$-polarized K3 surface
in the sense of \cite{MR544937}
is a pair $(\bfU,j)$
of a K3 surface $Y$
and a primitive lattice embedding
$
j \colon \bfU \hookrightarrow \Pic Y.
$
As explained, e.g., in \cite{Huybrechts_K3},
an elliptic K3 surface with a section
corresponds naturally to a pseudo-ample $\bfU$-polarized K3 surface.
Fix a primitive embedding of $\bfU$
to the K3 lattice
$
 \KL = \bfU \bot \bfU \bot \bfU \bot \bfE_8 \bot \bfE_8,
$
which is unique up to the left action of $\rO(\KL)$,
and let
$
 \bfT = \bfU \bot \bfU \bot \bfE_8 \bot \bfE_8
$
be the orthogonal lattice.
As explained in \cite[Section 3]{MR1420220},
the global Torelli theorem \cite{MR0284440, MR0447635}
and the surjectivity of the period map \cite{MR592693}
show that the period map gives an isomorphism
from the coarse moduli scheme of
pseudo-ample $\bfU$-polarized K3 surfaces
to the quotient
$
 M \coloneqq \Gamma \backslash \cD 
$
of the bounded Hermitian domain
\begin{align} \label{eq:D}
 \cD \coloneqq \lc \ld \Omega \rd \in \bfP(\bfT \otimes \bC)
  \relmid (\Omega, \Omega) = 0, \ (\Omega, \Omegabar) > 0 \rc
\end{align}
of type IV by
$
 \Gamma \coloneqq \rO(\bfT).
$

The moduli space of elliptic K3 surfaces with a section
attracts much attention recently,
not only from the point of view of modular compactification
(see e.g.~\cite{1902.10686,2002.07127}
and references therein),
but also because of the relation with tropical geometry
and mirror symmetry
\cite{MR3937674,1810.07685}.

A \emph{modular form}
on $\cD$
with respect to $\Gamma$
of weight $k \in \bZ$
and character $\chi \in \Char(\Gamma) \coloneqq \Hom(\Gamma, \bCx)$
is a holomorphic function
$
 f \colon \cDtilde \to \bC
$
on the total space
\begin{align}
  \cDtilde \coloneqq \lc \ld \Omega \rd \in \bfT \otimes \bC
  \relmid (\Omega, \Omega) = 0, \ (\Omega, \Omegabar) > 0 \rc
\end{align}
of a principal $\bCx$-bundle on $\cD$
satisfying
\begin{enumerate}[(i)]
 \item
$f(\alpha z) = \alpha^{-k} f(z)$
for any $\alpha \in \bCx$, and
 \item
$f(\gamma z) = \chi(\gamma) f(z)$
for any $\gamma \in \Gamma$.
\end{enumerate}
%\begin{align} \label{eq:mod1}
% f(\lambda z) = \lambda^{-k} f(z)
%\end{align}
%for any $\lambda \in \bCx$ and
%\begin{align} \label{eq:mod2}
% f(\gamma z) = \chi(\gamma) f(z)
%\end{align}
%for any $\gamma \in \Gamma$.
The vector spaces
$A_k(\Gamma, \chi)$
of modular forms
constitute the ring
\begin{align}
 \Atilde(\Gamma)
  \coloneqq \bigoplus_{k=0}^\infty \bigoplus_{\chi \in \Char(\Gamma)}
      A_k(\Gamma, \chi)
\end{align}
of modular forms.
We also write the subring of modular forms without characters as
\begin{align}
 A(\Gamma)
  \coloneqq \bigoplus_{k=0}^\infty
      A_k(\Gamma).
\end{align}

Let
$
V \coloneqq \Spec \bC[x,w]
$
be a 2-dimensional affine space over $\bC$.
For $k \in \bN$,
we write the $k$-th symmetric product of $V$ as $\Sym^k V$.
The special linear group $\SL_2$ acts naturally
on
$
S \coloneqq \Sym^8 V \times \Sym^{12} V
$
considered as an affine variety,
whose coordinate ring will be denoted by
\begin{align}
  \bC[S] = \bC[u_{8,0}, u_{7,1}, \ldots, u_{0,8}, u_{12,0}, u_{11,1}, \ldots, u_{0,12}].
\end{align}
We let $\Gm$ act on $S$
in such a way that
$u_{i,j}$ has weight $(i+j)/2$.
This $\Gm$-action commutes with the $\SL_2$-action,
so that the invariant subring
$
\bC[S]^{\SL_2}
$
has an induced $\Gm$-action.

Building on \cite{MR615858},
it is shown in \cite[Theorem 7.9]{1810.07685}
that the period map induces an isomorphism
from $\Proj \bC[S]^{\SL_2}$
to the Satake--Baily--Borel compactification
of $\Gamma \backslash \cD$,
so that one has an isomorphism
\begin{align} \label{eq:OO}
  A(\Gamma)
  \cong
  \bC [S]^{\SL_2}
\end{align}
of graded rings.

The main result of this paper is the following:

\begin{theorem} \label{th:main}
  One has
  \begin{align} \label{eq:main}
    \Atilde(\Gamma) \cong 
    \left.
    \lb
    \bC [S]^{\SL_2}
    \rb
    \ld s_{132} \rd
    \middle/ \lb s_{132}^2 - \Delta_{264} \rb \right.
  \end{align}
  where $s_{132}$ is an element of weight $132$ and
  $
  \Delta_{264}
  \in
  \bC[S]^{\SL_2}
  $
  is an element of weight $264$.
\end{theorem}

The proof is based on the construction of an algebraic stack
which is isomorphic to the orbifold quotient $\ld \cD \middle/ \rO(T) \rd$
in codimension 1.
The same strategy has been used in
\cite{1406.0332}
and
\cite{2005.00231}
to determine the rings of modular forms with characters
for the lattices
$\bfU \bot \bfU \bot \bfE_8$
and
$\bfU \bot \bfU \bot \bfA_1 \bot \bfA_1$
respectively.

The modular form $s_{132}$ is constructed
in \cite[Lemma 5.1]{MR2357689}.
It can also be obtained
either as the quasi pull-back
\cite[Theorem 8.2]{MR3184170}
% \begin{align}
%   \left. \frac{\Phi_{12}(z)}{\prod_{r\in R_{-2}/\{\pm 1\}} (z,r)} \right| _\cD
% \end{align}
of the Borcherds form $\Phi_{12}$
associated with the even unimodular lattice
% $\II_{2,26}$
of signature $(2,26)$
\cite[Section 10, Example 2]{MR1323986},
% where $R_{-2}$ is the set of $(-2)$-vectors of
% $\bfE_8 \cong \bfT^\bot \subset \II_{2,26}$
or by applying \cite[Theorem 10.1]{MR1323986}
to the nearly holomorphic modular form
\begin{align}
  \frac{1728 E_4}{E_4^3 - E_6^2} = \frac{1}{q} + 264 + 8244 q + 139520 q^2 + \cdots
\end{align}
where
\begin{align}
  E_4
  &= 1 + 240 \sum_{n=1}^\infty \frac{n^3 q^n}{1-q^n}
  = 1 + 240 q + 2160 q^2 + \cdots, \\
  E_6
  &= 1 - 504 \sum_{n=1}^\infty \frac{n^5 q^n}{1-q^n}
  = 1 - 504 q - 6632 q^2 + \cdots.
  % \Delta
  % &= E_4^3 - E_6^2
  % = 1728 q - 41472 q^2 + \cdots.
\end{align}
In particular,
it is a cusp form
with character $\det$
admitting an infinite product expansion.
See also \cite[Section 5]{MR3911795}
and references therein
for the case of the even unimodular lattice
of signature (2,10).

Since $\SL_2$ is reductive,
the invariant ring
$\bC[S]^{\SL_2}$
is finitely generated,
and
there exists an algorithm
for computing a finite generating set
(see e.g.~\cite{MR2667486} and references therein).
The element $\Delta_{264}$ can also be computed algorithmically, and
it is an interesting problem to describe them explicitly.

% This paper is organized as follows:

\section{The coarse moduli space
of \texorpdfstring{$\bfU$}{U}-polarized K3 surfaces}
 \label{sc:coarse moduli}

As is well known
(cf.~e.g.~\cite[Section 4]{MR2732092}),
a $\bfU$-polarized K3 surface
admits a Weierstrass model of the form
\begin{align} \label{eq:weierstrass1}
 z^2 = y^3 + g_2(x,w;u) y + g_3(x,w;u)
\end{align}
in $\bfP(1,4,6,1)$,
where
\begin{align}
 g_2(x,w;u)
 &= \sum_{i=0}^8 u_{8-i,i} x^{8-i} w^i \\
 &= u_{8,0} x^8 + u_{7,1} x^7 u + \cdots + u_{0,8} w^8, \\
 g_3(x,w;u)
 &= \sum_{i=0}^{12} u_{12-i,i} x^{12-i} w^i \\
 &= u_{12,0} x^{12} + u_{11,1} x^{11} u + \cdots + u_{0,12} w^{12}
\end{align}
for
\begin{align}
 u
 = ((u_{8,0}, \ldots, u_{0,8}), (u_{12,0}, \ldots, u_{0,12}))
 \in S.
\end{align}
The hypersurface in $\bfP(1,4,6,1)$
defined by \pref{eq:weierstrass1}
has a singularity
worse than rational double points
on the fiber at $a \in \bfP^1$
if and only if $\ord_a(g_2) \ge 4$ and
$\ord_a(g_3) \ge 6$
(see e.g.~\cite[Proposition I\!I\!I.3.2]{MR1078016}).
Let
$
U
\subset 
S
$
be the open subscheme
parametrizing hypersurfaces
with at worst rational double points.

The parameter $u$
describing a given $\bfU$-polarized K3 surface
is unique up to the action of $\SL_2 \times \Gm$,
where
$\Gm$ acts on
$
\bfP(1,4,6,1) \times \Sym^8 V \times \Sym^{12} V
$
by
\begin{align}
  \Gm \ni \lambda \colon
  ((x,y,z,w),(u_{i,j})_{i,j})
  \mapsto
  (x,\lambda^{2} y,\lambda^{3} z,w), (\lambda^{(i+j)/2} u_{i,j})_{i,j})
\end{align}
rescaling the holomorphic volume form
\begin{align}
  \Omega= \Res \frac{w dx \wedge dy \wedge dz}{z^2 - y^3 - g_2(x,w;u) y - g_3(x,w;u)}
\end{align}
as
\begin{align} \label{eq:Cx-action}
  \Omega_{\lambda u}
  &= \Res \frac{w dx \wedge d(\lambda^2 y) \wedge d (\lambda^3 z)}{(\lambda^3 z)^2 - (\lambda^2 y)^3 - g_2(x,w;\lambda \cdot u) (\lambda^2 y) - g_3(x,w;\lambda \cdot u)}
  = \lambda^{-1} \Omega_u.
\end{align}
The categorical quotient
$
T
\coloneqq U / \SL_2
$
is the coarse moduli scheme of pairs
$(Y, \Omega)$
consisting of a $\bfU$-polarized K3 surface $Y$
and a holomorphic volume form $\Omega \in H^0(\omega_Y)$ on $Y$.
The fact that the codimension of
$
S \setminus U
$
is greater than 2
implies
an isomorphism
\begin{align} \label{eq:C[T]}
  \bC[S]^{\SL_2} \cong \bC[T]
\end{align}
of graded rings.
Since the character of
$\bC[S]$
as a $\SL_2 \times \Gm$-module
is given by
\begin{align}
  \prod_{i=0}^{8}
  \left( 1-q^{2i-8} t^4 \right)^{-1}
  \prod_{i=0}^{12}
  \left( 1-q^{2i-12} t^6 \right)^{-1},
\end{align}
the Hilbert series of the invariant ring
% $
% \bC[S]^{\SL_2}
% $
is given by
\begin{align}
  \sum_{i=0}^\infty \dim \lb \bC[S]^{\SL_2} \rb_i \, t^i
  =
  \Res_{q=0}
  \left(
    (q^{-1} - q)
    \prod_{i=0}^{8}
    \left( 1-q^{2i-8} t^4 \right)^{-1}
    \prod_{i=0}^{12}
    \left( 1-q^{2i-12} t^6 \right)^{-1}
  \right)
\end{align}
as explained, e.g., in \cite[Section 4.4]{MR2004218}.
It follows from the global Torelli theorem and the surjectivity of the period map
that the period map induces a ring isomorphism
\begin{align} \label{eq:A(Gamma)}
  A(\Gamma)
  \simto
  \bC[T],
\end{align}
which preserves the grading by \pref{eq:Cx-action}.
The isomorphism
\pref{eq:OO}
follows from \pref{eq:C[T]} and \pref{eq:A(Gamma)}.

\section{Modular forms with characters}

The coarse moduli space
$
M \coloneqq \Gamma \backslash \cD
$
of $\bfU$-polarized K3 surfaces
is an open subvariety
of its Satake--Baily--Borel compactification
$
 \Proj A(\Gamma)
  \cong \bfP(4^9,6^{13}) \GIT \SL_2.
%  \cong \bfP(2,3,4,5,6).
$
Although $M$ and
the orbifold quotient
$
 \bM \coloneqq \ld \Gamma \backslash \cD \rd
$
are closely related,
the canonical morphism $\bM \to M$
is not an isomorphism
even in codimension 1.
In order to obtain an orbifold
which is isomorphic to $\bM$ in codimension 1
(so that
%the Picard groups and 
the total coordinate rings are isomorphic),
consider the stack
\begin{align}
  \bP \coloneqq \ld \bP \lb 4^{9},6^{13} \rb / \SL_2 \rd,
\end{align}
defined as the quotient of $\bC^{24} \setminus \bszero$
by the action of $\SL_2 \times \Gm$.
The morphism
$
 \bM \to M
$
lifts to a morphism
$
 \bM \to \bP,
$
which is an isomorphism in codimension 0,
since the generic stabilizers are $\lc \pm \id \rc$
on both sides.

Stabilizers of $\bM$ along divisors come from
reflections.
One divisor with a generic stabilizer
comes from the reflection with respect to a $(-2)$-vector
whose reflection hyperplane
corresponds to the locus
where the Picard lattice contains
$\bfU \bot \bfA_1$.
In order to describe this locus,
first consider the discriminant
\begin{align}
 h(x,w;u)
 \coloneqq 4 g_2(x,w;u)^3 + 27 g_3(x,w;u)^2
\end{align}
of
$
 y^3 + g_2(x,w;u) y + g_3(x,w;u)
$
as a polynomial in $y$,
which is homogeneous
of degree $24$ in $(x,w)$
and degree $12$ in $u$.
Note that the discriminant of a polynomial
$
 \sum_{i=0}^n a_i x^i w^{n-i}
$
with respect to $(x,w)$
is homogeneous of degree $2(n-1)$
in $\bZ[a_0, \ldots, a_n]$
if $\deg a_0 = \cdots = \deg a_n = 1$.
It follows that the discriminant
$k_{552}(u)$ of $h(x,w;u)$ with respect to $(x,w)$
is a homogeneous polynomial
of degree $2 \cdot 23 \cdot 12 = 552$ in $u$.
A general point
on the divisor $\bD_{552}$ of $\bP$
defined by $k_{552}(u)$
corresponds to the locus
where two fibers of Kodaira type $\mathrm{I}_1$ collapse
into one fiber.
This divisor has two components;
a general point on one component corresponds to the case
when there exists a point $p = [x:w] $ on $\bfP^1$ such that
neither $g_2$ nor $g_3$ vanishes at $p$,
and a general point on the other component
corresponds to the case
when both $g_2$ and $g_3$ vanishes at $p$.
In the former case,
the resulting singular fiber is of Kodaira type $\mathrm{I}_2$,
and the surface acquires an $A_1$-singularity.
In the latter case,
the resulting singular fiber is of Kodaira type I\!I,
and the surface does not acquire any new singularity.
%The defining equation of $D_1$ is $\Delta_T$.
The defining equation of the latter component
is the resultant of $g_2$ and $g_3$.
It is given as the determinant
\begin{align}
 r_{96}(t)=
\begin{vmatrix}
  u_{8,0} & u_{7,1} & \cdots & u_{0,8}\\
  & u_{8,0} & \cdots & u_{1,7} & u_{0,8} \\
  && \ddots && \ddots & \ddots \\
  &&& \ddots && \ddots & \ddots \\
  &&&& u_{8,0} & \cdots & u_{1,7} & u_{0,8} \\
  u_{12,0} & u_{11,1} &\cdots & \cdots & u_{1,11} & u_{0,12} \\
  % & u_{12,0} & \cdots & u_{2,10} & u_{1,11} & u_{0,12} \\
  & \ddots &&&& \ddots & \ddots \\
  && u_{12,0} & u_{11,1} & \cdots & \cdots & u_{1,11} & u_{0,12}
\end{vmatrix}
\end{align}
of the Sylvester matrix,
which is homogeneous of degree
\begin{align}
  12 \times 4 + 8 \times 6 = 96.
\end{align}
% (note that the product of diagonal entries is given by
% $
% u_{8,0}^{12} u_{12}^8
% $).
As shown in \cite[Lemma 6.1]{1406.0332},
the polynomial $k_{552}(t)$ is divisible by $r_{96}(t)^3$, and
the quotient
\begin{align} \label{eq:Delta_60}
 \Delta_{264}(t) \coloneqq k_{552}(t)/r_{96}(t)^3
\end{align}
defines the reflection hyperplane
along a $(-2)$-vector.

Recall from
\cite{MR2450211,MR2306040}
that the \emph{root construction}
is an operation
which adds a stabilizer
along a divisor.
Let $\bT$ be the stack
obtained from $\bP$
by the root construction of order 2
along the divisor on $\bP$
defined by $\Delta_{264}(t)$,
which is the quotient of the double cover of $\bP$
branched along $\Delta_{264}(t)$
by the group $G$ of deck transformations.
The Picard group of $\bT$
(or the $G$-equivariant Picard group of $\bP$)
is generated by the pull-back
$
 \cO_{\bT}(1) \coloneqq p^* \cO_{\bP}(1)
$
of the generator $\cO_{\bP}(1)$ of the Picard group of $\bP$
by the structure morphism
$
 p \colon \bT \to \bP
$
and the line bundle $\cO_{\bT}(\bD_{132})$
such that the space
$
 H^0(\cO_{\bT}(\bD_{132}))
$
is generated by an element $s_{132}$
satisfying
$
 s_{132}^2
  = \Delta_{264}
  \in H^0(\cO_{\bT}(264))
  \cong H^0(\cO_{\bP}(264)).
$
Note also that
$
 \omega_{\bP} \cong \bO_{\bP}(a)
$
where
\begin{align}
  a
  = - \sum_{i=0}^8 \deg u_{8-i,i} - \sum_{i=0}^{12} \deg u_{12-i,i}
  = - 9 \times 4 - 13 \times 6
  = - 114.
\end{align}
The ramification formula for the canonical bundle gives
\begin{align}
 \omega_{\bT}
  &\cong p^* \omega_{\bP} \otimes \cO_{\bT}(\bD_{132}) \\
  &\cong \cO_{\bT}(-114) \otimes \cO_{\bT}(132+(-132+\bD_{132})) \\
  &\cong \cO_{\bT}(18) \otimes \cO_{\bT}(-132+\bD_{132}).
   \label{eq:omega_T}
\end{align}
Note that $\cO_{\bT}(-132+\bD_{132})$ is an element
of order two in $\Pic \bT$.
By comparing \pref{eq:omega_T} with
\begin{align} \label{eq:omega}
 \omega_{\bM} \cong \cO_{\bM}( \dim \bM) \otimes \det = \cO_\bM(18) \otimes \det
\end{align}
which follows from (the proof of) \cite[Proposition 5.1]{1406.0332},
one concludes that $\bM$ has no further stabilizer along a divisor,
so that the lift $\bM \to \bT$ of $\bM \to \bP$ is an isomorphism
in codimension 1.
It follows that
the injective map
$
\bZ \times \Char(\Gamma)
\to \Pic \bM,
$
$
(i,\chi)
\mapsto
\cO_\bM(i) \otimes \chi
$
is surjective,
and the total coordinate ring
(also known as the Cox ring) of $\bM$ is given by
\begin{align} \label{eq:TCR}
 \bigoplus_{\cL \in \Pic \bM} H^0(\cL)
  \cong \bigoplus_{i=0}^\infty H^0(\cO_\bM(i))
   / (s_{132}^2-\Delta_{264}(t)).
\end{align}

\bibliographystyle{amsalpha}
\bibliography{bibs}

\end{document}